\documentclass[1p]{elsarticle}
\usepackage{amssymb}
\usepackage{amsfonts}

\newtheorem{theorem}{Theorem}
\newtheorem{conjecture}[theorem]{Conjecture}

\newtheorem{lemma}[theorem]{Lemma}

\newproof{pf}{Proof}

\begin{document}
\title{Distant total sum distinguishing index of graphs}

\author{Jakub Przyby{\l}o\fnref{fn1,fn2}}
\ead{jakubprz@agh.edu.pl, phone: 048-12-617-46-38,  fax: 048-12-617-31-65}

\fntext[fn1]{Financed within the program of the Polish Minister of Science and Higher Education
named ``Iuventus Plus'' in years 2015-2017, project no. IP2014 038873.}
\fntext[fn2]{Partly supported by the Polish Ministry of Science and Higher Education.}

\address{AGH University of Science and Technology, al. A. Mickiewicza 30, 30-059 Krakow, Poland}

\begin{abstract}
Let $c:V\cup E\to\{1,2,\ldots,k\}$ be a proper total colouring of
a graph $G=(V,E)$ with maximum degree $\Delta$.
We say vertices $u,v\in V$ are \emph{sum distinguished} if 
$c(u)+\sum_{e\ni u}c(e)\neq c(v)+\sum_{e\ni v}c(e)$.
By $\chi''_{\Sigma,r}(G)$ we denote the
least integer $k$ 
admitting such a colouring $c$
for which
every $u,v\in V$, $u\neq v$, at distance at most $r$ from each other are sum distinguished in $G$.
For every positive integer $r$ 
an infinite family of examples is known with $\chi''_{\Sigma,r}(G)=\Omega(\Delta^{r-1})$.
In this paper we prove that $\chi''_{\Sigma,r}(G)\leq (2+o(1))\Delta^{r-1}$ for every integer $r\geq 3$
and each graph $G$, while $\chi''_{\Sigma,2}(G)\leq (18+o(1))\Delta$.
\end{abstract}

\begin{keyword}
$r$-distant total neighbour sum distinguishing index of a graph
\end{keyword}

\maketitle

\section{Introduction}
In the paper~\cite{ChartrandErdosOellermann} Chartrand, Erd\H{o}s and Oellermann, proposed
definitions of a few graph families that might serve a role of \emph{irregular graphs}.
This research originated from the basic phenomenon that in fact there are no non-trivial irregular graphs at all,
which fulfill the most natural condition that all their vertex degrees are pairwise distinct.
Facing no definite obvious solution to this issue,
Chartrand et al.~\cite{Chartrand} turned towards measuring the ``level of irregularity'' of a graph $G=(V,E)$ via the following graph invariant.
The \emph{irregularity strength} of $G$ is the least $k$ so that we are able to construct an irregular multigraph (a multigraph with pairwise distinct vertex degrees) of $G$
by multiplying some of its edges - each at most $k$ times.
This parameter, denoted by $s(G)$, is well defined whenever $G$ contains at most one isolated vertex and no isolated edges.
Note that equivalently, 
$s(G)$ may be defined as the least $k$ so that an edge colouring $c:E\to\{1,2,\ldots,k\}$ exists with
$\sum_{w\in N(u)}c(uw)\neq \sum_{w\in N(v)}c(vw)$
for every $u,v\in V$, $u\neq v$.
This graph invariant has been analyzed in multiple papers, see e.g.~\cite{Aigner,Bohman_Kravitz,Lazebnik,
Faudree2,Faudree,Frieze,KalKarPf,Lehel,MajerskiPrzybylo2,Nierhoff,Przybylo},
but also gave rise to a fast-developing branch of research, which may be referred to as ``additive graph labellings'', or more generally ``vertex distinguishing colourings'', see~\cite{Aldred,Aigner3,ref_BacJenMilRya1,BurrisSchelp,Observability,FlandrinMPSW,123KLT,PilsniakWozniak_total,12Conjecture,Zhang_total,Zhang} for examples of papers introducing a few representative concepts of this branch. 

One of these is a problem of a \emph{total neighbour sum distinguishing colouring} of a given graph graph $G=(V,E)$,
i.e. a proper total colouring $c:E\to\{1,2,\ldots,k\}$ such that $c(u)+\sum_{w\in N(u)}c(uw)\neq c(v)+\sum_{w\in N(v)}c(vw)$
for every $uv\in E$, see~\cite{PilsniakWozniak_total}. The least $k$ admitting such a colouring $c$ is called the \emph{total neighbour sum distinguishing index} of $G$ and denoted by $\chi''_{\Sigma}(G)$. We also denote by 
$$w_c(v):=c(v)+\sum_{w\in N(v)}c(vw)$$ the so called \emph{weighted (total) degree} of a vertex $v\in V$ (note that $w_c$ is required to define a proper vertex colouring of $G$), and call $u,v\in V$ \emph{sum-distinguished} if $w_c(u)\neq w_c(v)$. 
It was conjectured in~\cite{PilsniakWozniak_total} that $\chi''_{\Sigma}(G)\leq \Delta(G)+3$ for every graph $G$ (note that $\chi''_{\Sigma}(G)\geq \Delta(G)+1$ due to the properness of the total colourings investigated). So far the best general upper bound implies that $\chi''_{\Sigma}(G)\leq \Delta(1+o(1))$ for every graph $G$ with maximum degree $\Delta$, see~\cite{Przybylo_asym_optim_total} and~\cite{Przybylo_asymptotic_note}. See also~\cite{total_sum_planar,LiLiuWang_sum_total_K_4,PilsniakWozniak_total,Przybylo_CN_3} for partial results concerning this conjecture.

In this paper we investigate a generalization of the problem above 
stemming among others from
the concept of
distant chromatic numbers, see e.g.~\cite{DistChrSurvey} for a survey concerning these.
Let $r$ be a positive integer. Vertices $u,v\in V$ shall be called $r$-neighbours if $1\leq d(u,v)\leq r$, where $d(u,v)$ denotes the distance of $u$ and $v$ in $G$. We say a proper total colouring $c:V\cup E\to\{1,2,\ldots,k\}$ is \emph{$r$-distant sum distinguishing} if $w_c(u)\neq w_c(v)$ whenever $u,v$ are $r$-neighbours in $G$. The least $k$ admitting such a colouring is denoted by $\chi''_{\Sigma,r}(G)$ and called the \emph{$r$-distant total sum distinguishing index} of $G$.
In~\cite{Przybylo_distant} 
a similar graph invariant, denoted by ${\rm ts}_r(G)$ was investigated (in fact it is just an analogous as above distant generalization of a total version of the irregularity strength of graphs introduced in~\cite{ref_BacJenMilRya1}). The only difference was that the total colourings investigated there did not have to be proper, hence obviously ${\rm ts}_r(G)\leq \chi''_{\Sigma,r}(G)$ for every graph $G$.
From this paper it follows that the general upper bound for ${\rm ts}_r(G)$, hence also for $\chi''_{\Sigma,r}(G)$ cannot be much smaller than $\Delta^{r-1}$, see~\cite{Przybylo_distant} for details.
In this paper on the other hand we prove that $\Delta^{r-1}$ is the right order of a general upper bound for $\chi''_{\Sigma,r}$ by proving Theorem~\ref{Main_Theorem_chi''Sigmar_probabil} below.


\section{Main Result and Probabilistic Preliminaries}
Our main result is the following.
\begin{theorem}\label{Main_Theorem_chi''Sigmar_probabil}
For every integer $r\geq 2$
there exists a constant $\Delta_0$ such that for each graph $G$ with maximum degree $\Delta\geq\Delta_0$,
\begin{equation}\label{main_total_distant_bound}
\chi''_{\Sigma,r}(G)\leq 2\Delta^{r-1}+5\Delta^{r-\frac{4}{3}}\ln^2\Delta+16\Delta+6,
\end{equation}
hence for all graphs
$$\chi''_{\Sigma,r}(G)\leq (2+o(1))\Delta^{r-1}$$
if $r\geq 3$, while for $r=2$:
$$\chi''_{\Sigma,2}(G)\leq (18+o(1))\Delta.$$
\end{theorem}
We note that the best thus far general upper bound for ${\rm ts}_r(G)$, a parameter upper-bounded by $\chi''_{\Sigma,r}(G)$, is also of the form $(2+o(1))\Delta^{r-1}$ (for $r\geq 2$), see~\cite{Przybylo_distant_total_irreg_probabil}.
The proof of Theorem~\ref{Main_Theorem_chi''Sigmar_probabil} 
is based on a colouring algorithm described in Section~\ref{section_with_algorithm}.
A different type of such an algorithm, but with a similar flavour
can in fact be found in~\cite{Przybylo_distant_total_irreg_probabil} 
(and some other more distant protoplasts e.g. in~\cite{AnhKalPrz,Kalkowski12,MajerskiPrzybylo1}).
Actually, we shall 
make use of a specific ordering of the vertices essentially proved to exist there, 
see Lemma~\ref{MainSequencingLemma_total} below. We include a proof of its existence though for the sake of completeness of the presentation of our reasoning.

We shall use a probabilistic approach based on
the Lov\'asz Local Lemma, see e.g.~\cite{AlonSpencer},
and the Chernoff Bound, see e.g.~\cite{JansonLuczakRucinski}
(Th. 2.1, page 26). 
\begin{theorem}[\textbf{The Local Lemma}]
\label{LLL-symmetric}
Let $A_1,A_2,\ldots,A_n$ be events in an arbitrary pro\-ba\-bi\-li\-ty space.
Suppose that each event $A_i$ is mutually independent of a set of all the other
events $A_j$ but at most $D$, and that ${\rm \emph{\textbf{Pr}}}(A_i)\leq p$ for all $1\leq i \leq n$. If
$$ ep(D+1) \leq 1,$$
then $ {\rm \emph{\textbf{Pr}}}\left(\bigcap_{i=1}^n\overline{A_i}\right)>0$.
\end{theorem}
\begin{theorem}[\textbf{Chernoff Bound}]\label{ChernofBoundTh}
For any $0\leq t\leq np$,
$${\rm\emph{\textbf{Pr}}}({\rm BIN}(n,p)>np+t)<e^{-\frac{t^2}{3np}}~~{and}~~{\rm\emph{\textbf{Pr}}}({\rm BIN}(n,p)<np-t)<e^{-\frac{t^2}{2np}}\leq e^{-\frac{t^2}{3np}}$$
where ${\rm BIN}(n,p)$ is the sum of $n$ independent Bernoulli variables, each equal to $1$ with probability $p$ and $0$ otherwise.
\end{theorem}
Note that if $n\leq k$, then we may still apply the Chernoff Bound above
to prove that
$\mathbf{Pr}(X>kp+t) < e^{-\frac{t^2}{3kp}}$ (for $t\leq\lfloor k\rfloor p$).


Given any graph $G=(V,E)$ of maximum degree $\Delta$, positive integer $r$ and a vertex $v\in V$, we denote by $N^r(v)$ the set of all $r$-neigbours of $v$ in $G$, and set $d^r(v)=|N^r(v)|$. We also define a partition of $V$ into:
%
\begin{eqnarray}
S&=& S(G) = \left\{u\in V:d(u)\leq \Delta^{\frac{2}{3}}\right\};\label{S_definition}\\
B&=& B(G) = \left\{u\in V:d(u)> \Delta^{\frac{2}{3}}\right\};\label{B_definition}
\end{eqnarray}
and denote $S(v)=\{u\in N(v):u\in S\}$, $s(v)=|S(v)|$, $B(v)=\{u\in N(v):u\in B\}$, $b(v)=|B(v)|$.
%
Moreover, given any fixed linear ordering of the vertices of $G$, every
neighbour or $r$-neighbour $u$ which precedes $v$ in 
this
ordering 
shall be called a \emph{backward neighbour} or \emph{$r$-neighbour}, resp., of $v$.
Analogously, the remaining ones shall be called \emph{forward neighbours} or \emph{$r$-neighbours}, resp., of $v$, while the edges joining $v$ with its forward or backward neighbours shall be referred to as \emph{forward} or \emph{backward}, respectively.
Additionally, for any subset $U\subset V$, let $N_-(v)$, $N_-^r(v)$, $N_U^r(v)$ denote the sets of all backward neighbours, backward $r$-neighbours and $r$-neighbours in $U$ of $v$, respectively.
Set $d_-^r(v)=|N_-^r(v)|$, $d_U^r(v)=|N_U^r(v)|$, and let $b_-(v)$ denote the number of backward neighbours of $v$ which belong to $B(v)$.
Note also that if we denote $D(v)=\sum_{u\in N(v)}d(u)$, then for $r\geq 2$:
$$d^r(v) \leq D(v)\Delta^{r-2} \leq d(v)\Delta^{r-1}.$$

\begin{lemma}
\label{MainSequencingLemma_total}
For every integer $r\geq 2$ there exists $\Delta'_0$ such that given any graph $G=(V,E)$ with maximum degree $\Delta\geq \Delta'_0$, we may assign to every vertex $v\in V$ a real number $x_v\in[0,1]$ so that for any ordering $v_1,v_2,\ldots,v_n$ of the vertices of $G$ such that $x_{v_i}\leq x_{v_j}$ whenever $i\leq j$, 
it holds that for every vertex $v$ in $G$ with $b(v)\geq \Delta^{\frac{1}{3}}\ln\Delta$:
\begin{itemize}
\item[(i)] $d^r_I(v) \leq 2d(v)\Delta^{r-\frac{4}{3}} \ln\Delta$;
\item[(ii)] if $v\in R$, then: $b_-(v)\geq x_v b(v)-\sqrt{x_v b(v)}\ln\Delta$;
\item[(iii)] if $v\in R$, then: $d^r_-(v)\leq x_v D(v) \Delta^{r-2}+\sqrt{x_vD(v)\Delta^{r-2}}\ln\Delta$
\end{itemize}
where
$$I=\left\{v\in V: x_v<\frac{\ln\Delta}{\Delta^{\frac{1}{3}}}\right\};~~~~~~
R=\left\{v\in V : x_v\geq\frac{\ln\Delta}{\Delta^{\frac{1}{3}}}\right\}.$$
\end{lemma}


\begin{pf}
Fix $r\geq 2$ and $G=(V,E)$.
We do not specify $\Delta'_0$, but assume that $\Delta$ 
is large enough so that all explicit inequalities below hold.
We order the vertices of $V$ randomly. 
For this goal, we
associate with every vertex $v\in V$ an (independent) random variable $X_v\sim U[0,1]$ having the uniform distribution on $[0,1]$ (whose random value $x_v$ determines whether $v\in R$ or not).
We then consider any ordering $v_1,v_2,\ldots,v_n$ of the vertices in $V$ such that
$X_{v_i}\leq X_{v_j}$ ($x_{v_i}\leq x_{v_j}$) whenever $i\leq j$, $i,j\in\{1,2,\ldots,n\}$.

Let $v$ be a vertex of degree $d$ in $G$ with $b(v)\geq \Delta^{\frac{1}{3}}\ln\Delta$ (thus $d\geq \Delta^{\frac{1}{3}}\ln\Delta$).
Let then $A_{v,1}$ denote the event that $d^r_I(v) > 2d\Delta^{r-\frac{4}{3}} \ln\Delta$,
let $A_{v,2}$ be the event that $v$ belongs to $R$ and $b_-(v) < X_v b(v)-\sqrt{X_v b(v)}\ln\Delta$, and
let $A_{v,3}$ denote the event that $v$ belongs to $R$ and $d^r_-(v) > X_v D(v) \Delta^{r-2}+\sqrt{X_vD(v)\Delta^{r-2}}\ln\Delta$.

Since for every $u\in N^r(v)$,
the probability that $u$ belongs to $I$ equals $\frac{\ln\Delta}{\Delta^{\frac{1}{3}}}$,
and $|N^r(v)| \leq d\Delta^{r-1}$,
by the Chernoff Bound, 
\begin{eqnarray}
\mathbf{Pr}(A_{v,1}) &\leq& \mathbf{Pr}\left(d^r_I(v) > d\Delta^{r-\frac{4}{3}} \ln\Delta + \sqrt{d\Delta^{r-\frac{4}{3}} \ln\Delta}\ln\Delta\right) \nonumber\\
&<& e^{-\frac{d\Delta^{r-\frac{4}{3}} \ln^3\Delta}{3d\Delta^{r-\frac{4}{3}} \ln\Delta}} = \Delta^{-\frac{\ln\Delta}{3}} < \frac{1}{\Delta^{3r}}.\label{Av1Ineq_total}
\end{eqnarray}
Note also that by the Chernoff Bound, for any $x\in [0,1]$ for which $\sqrt{xb(v)}\ln\Delta \leq x b(v)$,
\begin{eqnarray}
& & \mathbf{Pr}(b_-(v) < X_v b(v)-\sqrt{X_v b(v)}\ln\Delta | X_v=x)\nonumber\\
&=& \mathbf{Pr}({\rm BIN}(b(v),x) < x b(v)-\sqrt{x b(v)}\ln\Delta) < \frac{1}{\Delta^{3r}}.\nonumber
\end{eqnarray}
As trivially $\mathbf{Pr}(b_-(v) < X_v b(v)-\sqrt{X_v b(v)}\ln\Delta | X_v=x)=0$ for any $x\in [0,1]$ such that $\sqrt{xb(v)}\ln\Delta > x b(v)$, we conclude that:
\begin{equation}
\mathbf{Pr}(A_{v,2}) \leq \mathbf{Pr}(b_-(v) < X_vb(v)-\sqrt{X_v b(v)}\ln\Delta) \leq \int\limits_{0}^1\frac{1}{\Delta^{3r}}dx=\frac{1}{\Delta^{3r}}.\label{Av2Ineq_total}\end{equation}
Finally
note 
that for $x\in[0, \frac{\ln\Delta}{\Delta^{\frac{1}{3}}})$,
\begin{equation}
\mathbf{Pr}(d^r_-(v) > X_v D(v) \Delta^{r-2}+\sqrt{X_vD(v)\Delta^{r-2}}\ln\Delta \wedge v\in R | X_v=x) = 0, \label{Av3Ineq_total_Part1}
\end{equation}
while
for $x\in [\frac{\ln\Delta}{\Delta^{\frac{1}{3}}},1]$,
analogously as above, by the Chernoff Bound:
\begin{eqnarray}
&&\mathbf{Pr}(d^r_-(v) > X_v D(v) \Delta^{r-2}+\sqrt{X_vD(v)\Delta^{r-2}}\ln\Delta \wedge v\in R | X_v=x) \nonumber\\
&\leq& \mathbf{Pr}({\rm BIN}(D(v) \Delta^{r-2},x) > x D(v) \Delta^{r-2}+\sqrt{xD(v)\Delta^{r-2}}\ln\Delta) <\frac{1}{\Delta^{3r}}\label{Av3Ineq_total_Part2}
\end{eqnarray}
(as $x\geq \frac{\ln\Delta}{\Delta^{\frac{1}{3}}}$ and $b(v)\geq \Delta^{\frac{1}{3}}\ln\Delta$, where $D(v)\geq b(v)\Delta^{\frac{2}{3}}$, imply that $\sqrt{xD(v)\Delta^{r-2}}\ln\Delta \leq x D(v) \Delta^{r-2}$).
By~(\ref{Av3Ineq_total_Part1}) and~(\ref{Av3Ineq_total_Part2}),
\begin{equation}
\mathbf{Pr}(A_{v,3})
\leq \int\limits_{0}^1\frac{1}{\Delta^{3r}}dx=\frac{1}{\Delta^{3r}}.\label{Av3Ineq_total}\end{equation}

Observe now that each event $A_{v,i}$ is mutually independent of all other events except those $A_{u,j}$ with $u$ at distance at most $2r$ from $v$, $i,j\in\{1,2,3\}$,
i.e., except at most $3\Delta^{2r}+2$ other events.
On the other hand, by~(\ref{Av1Ineq_total}), (\ref{Av2Ineq_total}) and~(\ref{Av3Ineq_total}), we have that $\mathbf{Pr}(A_{v,i}) < \Delta^{-3r}$ for every $v\in V$ with $b(v)\geq \Delta^{\frac{1}{3}}\ln\Delta$ and each $i\in\{1,2,3\}$.
Thus by the Lov\'asz Local Lemma, with positive probability none of such events $A_{v,i}$ 
appears, hence there exist choices of $x_v$ for all $v\in V$ for which the thesis holds.
\qed
\end{pf}

\section{Proof of Theorem~\ref{Main_Theorem_chi''Sigmar_probabil}}\label{section_with_algorithm}
Fix an integer $r\geq 2$. Within our proof we shall not specify $\Delta_0$.
Instead, we shall again assume that $G=(V,E)$ is a graph with sufficiently large maximum degree $\Delta$,
i.e. large enough so that all inequalities below are fulfilled, and not smaller than $\Delta'_0$ from Lemma~\ref{MainSequencingLemma_total} above.

Let $v_1,v_2,\ldots,v_n$ be an ordering of the vertices of $V$ guaranteed by Lemma~\ref{MainSequencingLemma_total},
with $x_v$ for $v\in V$,
$I$ and $R$ as specified in this lemma (and let $S,B$ be the vertex subsets defined in~(\ref{S_definition}) and~(\ref{B_definition}) in $G$).
Set
$$k:=\lceil\Delta^{r-\frac{4}{3}}\ln^2\Delta\rceil ~~~~{\rm and}~~~~ K:=\left\lceil\frac{\Delta^{r-1} +6\Delta +k}{k}\right\rceil \cdot k$$ 
and note that $k$ divides $K$ and that $\Delta^{r-1} +6\Delta +k\leq K\leq \Delta^{r-1} +6\Delta +2k$.
Let 
\begin{eqnarray}
L&=&\left(\bigcup_{l=1}^{\lceil\frac{\Delta+1}{k}\rceil-1}\left[K+(l-1)4k+1,K+(l-1)4k+k\right]\right)\nonumber\\
&\cup& \left[K+\left(\left\lceil\frac{\Delta+1}{k}\right\rceil-1\right)4k+1,K+\left(\left\lceil\frac{\Delta+1}{k}\right\rceil-1\right)4k
+\left(\Delta+1\right)-\left(\left\lceil\frac{\Delta+1}{k}\right\rceil-1\right)k\right].\nonumber
\end{eqnarray}
(In fact $L$ consists only of the last interval for every $r\geq 3$.)
Note that $L\subseteq[K+1,K+4\Delta+1]$, $L$ contains exactly $\Delta+1$ integers and for any two integers $i_1,i_2\in L$:
\begin{equation}\label{i1i2distinction}
i_1\neq i_2 \Rightarrow \left(\{i_1-k,i_1,i_1+k,i_1+2k\}\cap \{i_2-k,i_2,i_2+k,i_2+2k\} = \emptyset \pmod K\right),
\end{equation}
i.e. for $i_1\neq i_2$ no element of the four element set corresponding to $i_1$ is congruent modulo $K$ to any element from the set associated with $i_2$.
We first apply Vizing's theorem to colour properly the edges of $G$ with integers from $L$.
Then we greedily colour the vertices of $G$ with colours $1,2,\ldots,K$ so that the obtained total colouring of $G$ is proper modulo $K$
(i.e. the colours are not congruent modulo $K$ for the adjacent vertices, adjacent edges and for the edges and their ends).
Now we shall use a recoloring algorithm eventually yielding our final proper total colouring
$f:V\cup E\to\{1,2,\ldots,2K+k+4\Delta+1\}$.
Within this algorithm we shall be analyzing consecutive vertices in the fixed ordering (starting from $v_1$).
By $c_t(a)$ we shall denote a contemporary colour of every $a\in V\cup E$ within 
the ongoing algorithm.
After every stage of the algorithm (i.e. after analyzing a consecutive vertex in the sequence) corresponding to $v_i$,
the obtained total colouring of $G$ shall be required to remain proper modulo $K$,
except for some possible colour conflicts between vertices $v_j$ with $j>i$ and their adjacent vertices or incident edges,
which shall be taken care of later on. (In fact the admitted alterations of edge colours 
within the algorithm combined with property~(\ref{i1i2distinction}) shall immediately guarantee that our colouring restricted to edges shall always be proper modulo $K$, cf. ($1^\circ$)-($4^\circ$) below.)
The final target sum of every vertex $v\in V$, $w_f(v)$, shall be chosen and fixed the moment $v$
is analyzed so that $w_f(v)\neq w_f(u)$ for every backward $r$-neighbour $u$ of $v$ in $G$.
Ever since this moment we shall also require that after every succeeding stage, the total sum of $v$, $w_{c_t}(v)$, equals $w_f(v)$.
The moment a given vertex $v$ is analyzed we shall choose its (new) colour $c'(v)\in\{1,2,\ldots,K\}$ (and set $c_t(v)=c'(v)$ at this point).
Ever since this moment we shall admit four possible colours for $v$, namely from then on we shall always require so that
\begin{equation}\label{ctvset}
  c_t(v)\in\{c'(v),c'(v)+k,c'(v)+K,c'(v)+K+k\}
\end{equation}
(where $c'(v)+K+k\leq 2K+k$).
We shall admit at most two alterations of the colour for every edge in $E$ during the algorithm - only when any of its ends is being
analyzed.
Let $v$ be a currently analyzed vertex in the sequence $v_1,v_2,\ldots,v_n$, and let $u\in N(v)$ be any of its neighbours.
We admit the following alterations at this stage of our construction (associated with adjusting $w_f(v)$):
\begin{itemize}
\item[($1^\circ$)] adding $0$ or $K$ to $c_t(uv)$ if $uv$ is a forward edge of $v$, $v\in S$ and $u\in B$,
\item[($2^\circ$)] adding $0,$ or $k$ to $c_t(uv)$ if $uv$ is a forward edge of $v$ and $v\in B$ or $u\in S$,
\item[($3^\circ$)] adding $-K$, $0$ or $K$ to $c_t(uv)$ if $uv$ is a backward edge of $v$ and $u\in B$,
\item[($4^\circ$)] adding $-k$, $0$ or $k$ to $c_t(uv)$ if $uv$ is a backward edge of $v$ and $u\in S$,
\item[($5^\circ$)] adding $-K$, $0$ or $K$ to $c_t(u)$ if $u$ is a backward neighbour of $v$ and $u\in B$,
\item[($6^\circ$)] adding $-k$, $0$ or $k$ to $c_t(u)$ if $u$ is a backward neighbour of $v$ and $u\in S$,
\item[($7^\circ$)] adding any integer in $[-K+1,K-1]$ to $c_t(v)$.
\end{itemize}
These however have to performed so that all our requirements described above ar met at the end of this step;
in particular so that $c_t(u)\in\{c'(u),c'(u)+k,c'(u)+K,c'(u)+K+k\}$ and $w_{c_t}(u) = w_f(u)$ afterwards if $u$ is a backward neighbour of $v$,
hence we have always only two admissible additions via 
($5^\circ$) and ($6^\circ$) including $0$ due to the first of this conditions (cf.~(\ref{ctvset})), and
consequently also only two associated options via each of ($3^\circ$) and ($4^\circ$) so that also $w_{c_t}(u) = w_f(u)$ holds at the end of the step.
Note that the admitted alterations guarantee that $c_t(e)\in\{1,2,\ldots,2K+k+4\Delta+1\}$ for every $e\in E$ at every stage of the construction.

Suppose we are about to analyze a vertex $v=v_i$, $i\in\{1,2,\ldots,n\}$, and thus far all our requirements have been fulfilled, and all rules above have been obeyed.
We shall show that in every case the admitted alterations
provide us more options
(guaranteeing also that there shall be no colour conflicts between $v$ and its incident edges and adjacent vertices now and in the future provided that we follow the same rules and constraints in the further part of the algorithm)
for $w_{c_t}(v)$
than there are backward $r$-neighbours of $v$, and hence one of this options can be fixed as $w_f(v)$ so that this value is distinct from every $w_f(u)$ already fixed for any $u\in N^r_-(v)$. Denote the degree of $v$ by $d$,
and assume that $d>0$ (otherwise, we set $w_f(v)=1$ and $c_t(v)=1$).
Note that due to the admitted alterations, by the required properness modulo $K$ of the total colouring after every stage of the construction (with the exception for the forward vertices), applying ($7^\circ$) to $v$
we have to choose $c'(v)\in[1,K]$ so that neither of $c'(v)$ and $c'(v)+k$ is equal to $c'(u)$ or $c'(u)+k$ modulo $K$ for any backward $r$-neighbour $u$ of $v$ (this excludes at most $3\Delta$ of the $K$ potential options for $c'(v)$).
Analogously, we must choose $c'(v)\in[1,K]$ so that after this step neither of $c'(v)$ and $c'(v)+k$ equals the current and any of the future colours of the incident edges of $v$ modulo $K$. If $v\in B$ there are however only two distinct modulo $K$ colours available for every edge incident with $v$ now and in the further part of the construction; this yield additional $3\Delta$ constraints in this case.
If on the other hand $v\in S$, neither of $c'(v)$ and $c'(v)+k$ may be congruent to any element of the set
$c_t(e)-k, c_t(e), c_t(e)+k$ or $c_t(e)+2k$ modulo $K$ for any edge $e$ incident with $v$, 
where $c_t(e)$ is the contemporary colour of $e$ before this stage, what yields at most $5d(v)\leq5\Delta^{\frac{2}{3}}<3\Delta$ additional constraints.
Note that for any $q\geq 0$ we have $(\lfloor qk/K\rfloor+1) (K-6\Delta) \geq qk \frac{K-6\Delta}{K}\geq 0.1 qk \geq  0.1q\Delta^{r-\frac{4}{3}}\ln^2\Delta$,
and consider the following cases:

\textbf{Case 1:} If $v\in I$, $v\in B$ and $b(v)\geq \Delta^{\frac{1}{3}}\ln\Delta$,
then the admitted alterations provide at least
$(\lfloor dk/K\rceil+1) (K-6\Delta) \geq 
0.1d\Delta^{r-\frac{4}{3}}\ln^2\Delta$
available options for $w_{c_t}(v)$.
As by ($i$) (from Lemma~\ref{MainSequencingLemma_total}), $|N^r_-(v)|\leq d_I^r(v)\leq 2d\Delta^{r-\frac{4}{3}} \ln\Delta < 0.1d\Delta^{r-\frac{4}{3}}\ln^2\Delta$,
at least one of these available options is distinct from all $w_f(u)$ with $u\in N^r_-(v)$.


\textbf{Case 2:} If $v\in S$, then the admitted alterations provide at least
$(\lfloor (s(v)k+b(v) K)/K\rfloor+1) (K-6\Delta) = (\lfloor s(v)k/K \rfloor + 1) (K-6\Delta) + b(v) (K-6\Delta) \geq 0.1s(v)\Delta^{r-\frac{4}{3}}\ln^2\Delta + b(v)(\Delta^{r-1}+\Delta^{r-\frac{4}{3}}\ln^2\Delta)$
available options for $w_{c_t}(v)$.
On the other hand,
$|N^r_-(v)|\leq d^r(v) \leq D(v)\Delta^{r-2} \leq (s(v)\Delta^{\frac{2}{3}}+b(v)\Delta)\Delta^{r-2}$,
hence at least one of these available options is distinct from all $w_f(u)$ with $u\in N^r_-(v)$.


\textbf{Case 3:} If $v\in B$ and $b(v) < \Delta^{\frac{1}{3}}\ln\Delta$, then the admitted alterations provide at least
$(\lfloor dk/K\rfloor+1) (K-6\Delta)\geq 0.1 d\Delta^{r-\frac{4}{3}}\ln^2\Delta$ available options for $w_{c_t}(v)$.
On the other hand, analogously as in the case above,
$|N^r_-(v)|\leq d^r(v) \leq
s(v)\Delta^{\frac{2}{3}}\Delta^{r-2} + b(v)\Delta^{r-1} < d\Delta^{r-\frac{4}{3}}+\Delta^{r-\frac{2}{3}}\ln\Delta < 0.1d\Delta^{r-\frac{4}{3}}\ln^2\Delta$,
as $v\in B$ implies that $d\geq\Delta^{\frac{2}{3}}$.
We thus have at least one option available for $v$ distinct from all $w_f(u)$ with $u\in N^r_-(v)$.


\textbf{Case 4:} If $v\in R$, $v\in B$ and $b(v) \geq \Delta^{\frac{1}{3}}\ln\Delta$, then by ($ii$) the number of available options for $w_{c_t}(v)$ via admitted alterations of colours of the edges incident with $v$ is not smaller than:
\begin{eqnarray}
&    & (\lfloor(b_-(v)K + (d-b_-(v))k)/K\rfloor+1) (K-6\Delta) \nonumber\\
&  = & b_-(v) (K-6\Delta) + (\lfloor (d-b_-(v))k/K\rfloor+1) (K-6\Delta) \nonumber\\
&\geq& b_-(v)(\Delta^{r-1} + \Delta^{r-\frac{4}{3}}\ln^2\Delta) + 0.1(d - b_-(v)) \Delta^{r-\frac{4}{3}}\ln^2\Delta \nonumber\\
&\geq& b_-(v)\Delta^{r-1} + 0.1d \Delta^{r-\frac{4}{3}}\ln^2\Delta \nonumber\\
&\geq& (x_v b(v)-\sqrt{x_v b(v)}\ln\Delta)\Delta^{r-1} + 0.1 d \Delta^{r-\frac{4}{3}}\ln^2\Delta \nonumber\\
&\geq& x_v b(v)\Delta^{r-1} - \sqrt{d}\Delta^{r-1}\ln\Delta + 0.1 d \Delta^{r-\frac{4}{3}}\ln^2\Delta \nonumber\\
&\geq& x_v b(v)\Delta^{r-1} + 0.1 d \Delta^{r-\frac{4}{3}}\ln^2\Delta - d \Delta^{r-\frac{4}{3}}\ln\Delta \nonumber
\end{eqnarray}
(where the last inequality follows by the fact that $d\geq\Delta^{\frac{2}{3}}$).
This number is however greater than the number of backward $r$-neighbours of $v$, as by ($iii$),
\begin{eqnarray}
|N^r_-(v)| &\leq& x_v D(v) \Delta^{r-2}+\sqrt{x_vD(v)\Delta^{r-2}}\ln\Delta \nonumber\\
&\leq& x_v (b(v)\Delta +
s(v)\Delta^{\frac{2}{3}}) \Delta^{r-2}+\sqrt{d\Delta^{r-1}}\ln\Delta \nonumber\\
&\leq& x_v b(v)\Delta^{r-1} + d \Delta^{r-\frac{4}{3}} + d \Delta^{r-\frac{4}{3}}\ln\Delta. \nonumber
\end{eqnarray}
Thus in all cases there is at least one available sum, say $w^*$, for $v$ which is distinct from all $w_f(u)$ with $u\in N^r_-(v)$.
We then set $w_f(v)=w^*$ and perform the admitted alterations (following all the rules stated above) on the edges incident with $v$ and on the vertices in $\{v\}\cup N_-(v)$
so that $w_{c_t}(v)=w^*$ afterwards. 

By our construction, after analyzing $v_n$, all $w_f(v_i)$ are fixed for $i=1,\ldots,n$ so that $w_f(u)\neq w_f(v)$ whenever $u$ and $v$ are $r$-neighbours in $G$.
We then set  $f(a)=c_t(a)$ for every $a\in V\cup E$, completing the construction of a desired total colouring $f$ of $G$.
Note that $1\leq f(e)\leq 2K+k+4\Delta+1$ for every $e\in E$ and 
$1\leq f(v)\leq 2K+k$ for every $v\in V$, hence the thesis follows.
\qed

\section{Remarks}

The character of the problem investigated changes for $r\geq 2$ (compared to $r=1$).
Within our approach 
we have put an effort to optimize the second order term (for $r\geq 3$) in the upper bound~(\ref{main_total_distant_bound}),
up to a constant and the power in the logarithmic factor.
The result for $r=2$ could still be improved within the same approach, but we focused on the most general result with respect to all $r\geq 2$.
%
The case of $r=2$ involves special difficulties, what is reflected in our conjectures below, the first of which seems to be most challenging one in this field, while the second one is just its relaxation for the case of $r\geq 3$.
%
\begin{conjecture}\label{przybylo_main_con_sum_total_distant}
For every positive integer $r$ there exists a constant $C$ such that for each
graph $G$ with maximum degree $\Delta$,
$$\chi''_{\Sigma,r}(G) \leq \Delta+\Delta^{r-1}+C.$$
\end{conjecture}
\begin{conjecture}\label{przybylo_main_con_sum_total_distant_2}
For every integer $r\geq 3$ and for each
graph $G$ with maximum degree $\Delta$,
$$\chi''_{\Sigma,r}(G) \leq (1+o(1))\Delta^{r-1}.$$
\end{conjecture}

\end{document}